\documentclass[a4paper,10pt]{article}
\usepackage[utf8x]{inputenc}
\usepackage{amsthm,amsmath,amssymb,oldgerm}
\usepackage[margin=3cm]{geometry}
\theoremstyle{plain}
\newtheorem{thm}{Theorem}[section]
\newtheorem{lem}[thm]{Lemma}
\newtheorem{prop}[thm]{Proposition}
\newtheorem{cor}[thm]{Corollary}
\theoremstyle{definition}
 
\theoremstyle{definition}
\newtheorem{rem}[thm]{Remark}

\let\Im\relax
\DeclareMathOperator{\Im}{Im}
\DeclareMathOperator{\del}{\Delta_{hyp}}
\DeclareMathOperator{\hyp}{\mu_{hyp}} 
\DeclareMathOperator{\can}{\mu_{can}}

\DeclareMathOperator{\hatcan}{\widehat{\mu}_{can}}
\DeclareMathOperator{\shyp}{\mu_{shyp}}

\DeclareMathOperator{\kh}{{\it{K_{\mathbb{H}}}}}
\DeclareMathOperator{\khyp}{{\it{K_{\mathrm{hyp}}}}}
\DeclareMathOperator{\hkhyp}{{\it{HK}_{\mathrm{hyp}}}}

\DeclareMathOperator{\ghyp}{\it{g_{\mathrm{hyp}}}}
\DeclareMathOperator{\gcan}{\it{g_{\mathrm{can}}}}  
\DeclareMathOperator{\gh}{\it{g_{\mathbb{H}}}}

\let\Re\relax
\DeclareMathOperator{\Re}{Re}
\DeclareMathOperator{\vx}{\mathrm{vol_{\mathrm{hyp}}}}
\let\id\relax
\DeclareMathOperator{\id}{\mathrm{id}}
\DeclareMathAlphabet{\mathpzc}{OT1}{pzc}{m}{it}
\makeatletter

\newcommand{\Rmnum}[1]{\expandafter\@slowromancap\romannumeral #1@}
\title{Bounds for canonical Green's function at cusps}
{\small\author{Anilatmaja Aryasomayajula}}
\date{}
\begin{document}
\maketitle
\begin{abstract}
\noindent 
In this article, we derive bounds for the canonical Green's function defined on a noncompact hyperbolic 
Riemann surface, when evaluated at two inequivalent cusps. 
 
\vspace{0.2cm}\noindent
\end{abstract}
\section*{Introduction}
\paragraph{Notation and Main results}
Let $X$ be a noncompact hyperbolic Riemann surface of finite volume $\vx(X)$ with genus $g\geq 1$. Then, 
from the uniformization theorem from complex analysis, $X$ can be realized as the quotient space 
$\Gamma\backslash\mathbb{H}$, where $\Gamma\subset \mathrm{PSL}_{2}(\mathbb{R})$ is a Fuchsian 
subgroup of the first kind acting on the hyperbolic upper half-plane $\mathbb{H}$, via fractional linear 
transformations. Let $\mathcal{P}$ denote the set of cusps of $\Gamma$. Put $\overline{X}=X\cup \mathcal{P}$.  

\vspace{0.2cm}
The compact Riemann surface $\overline{X}$ is embedded in its Jacobian variety $\mathrm{Jac}(\overline{X})$ via the 
Abel-Jacobi map. Then, the pull back of the flat Euclidean metric by the Abel-Jacobi map is called the 
canonical metric, and the (1,1)-form associated to it is denoted by $\hatcan(z)$. We denote its 
restriction to $X$ by $\can(z)$. Put 
\begin{align*}
 d_{X}=\sup_{z\in X}\frac{\can(z)}{\shyp(z)}. 
\end{align*}
The canonical Green's function is defined as the unique solution of the differential equation (which is to 
be interpreted in terms of currents)
\begin{align*}
d_{z}d^{c}_{z}\gcan(z,w)+ \delta_{w}(z)=\can(z),
\end{align*}
with the normalization condition
\begin{align*}
\int_{X}\gcan(z,w)\can(z)=0. 
\end{align*}
Let $c_{X}$ denote a certain constant related to the Selberg zeta function (see equation \eqref{selbergconstant} for 
definition). Let $\lambda_{1}$ denote the first non-zero eigenvalue of the hyperbolic Laplacian $\del$ acting 
on smooth functions defined on $X$. Let $\kappa_{p}(z)$ denote the Kronecker's limit function associated to a cusp $p\in \mathcal{P}$, which is 
the constant term in the Laurent expansion of the Eisenstein series associated the cusp $p\in \mathcal{P}$ 
at $s=1$. Let $k_{p,q}(0)$ $\in$ $\mathbb{C}$ denote the zeroth Fourier coefficient of $\kappa_{p}(z)$ at 
the cusp $q\in\mathcal{P}$. 

\vspace{0.2cm}
With notation as above, we have the following upper bound
\begin{align*}
 &\big{|} g_{\mathrm{can}}(p,q)\big{|}\leq 4\pi \big{|}k_{p,q}(0)\big{|}+
\frac{2\pi}{g}\bigg(\sum_{\substack{s\in \mathcal{P}\\s\not = p}}\big{|}k_{s,p}(0)\big{|}+
\sum_{\substack{s\in  \mathcal{P}\\s\not = q}}\big{|}k_{s,q}(0)\big{|}\bigg)+\notag\\&
\frac{1}{\vx(X)}\bigg(\frac{4\pi(d_{X}+1)^{2}}{\lambda_{1}}+\frac{\big{|}4\pi c_{X}\big{|}}{g}+
\frac{ 43|\mathcal{P}|}{g}+4\pi\bigg)+\frac{2\log (4\pi)}{g}.
\end{align*}
\paragraph{Arithmetic significance}
Bounds for the canonical Green's function are very essential for calculating various arithmetic 
invariants like the Faltings height function and the faltings delta function. Especially bounds for the canonical Green's 
function evaluated at two inequivalent cusps are essential for calculating the arithmetic self-intersection number of 
the dualizing sheaf defined on an arithmetic surface. 

\vspace{0.2cm}
In \cite{abbes}, while bounding the arithmetic self-intersection number of 
the dualizing sheaf defined on the modular curve $X_{0}(N)$, A. Abbes and E. Ullmo computed bounds for the 
canonical Green's function evaluated at the cusps $0$ and $\infty$. In \cite{hartiwg}, H.~Mayer has done the same for 
the modular curve $X_{1}(N)$. 

\vspace{0.2cm}
Furthermore,  in \cite{kuhn} U.~K\"uhn as also derived bounds for the arithmetic self-intersection number 
of the dualizing sheaf defined on any curve defined over a number field. Our bounds hold true for any 
noncompact hyperbolic Riemann surface of genus $g>0$. So they can be directly 
used in \cite{kuhn}, and we hope that this leads to better bounds for U.~k\"uhn.

\vspace{0.2cm}
Lastly, using results from \cite{anilpaper1} and \cite{anilpaper2}, our bounds can be easily extended to 
the case when $X$ admits elliptic fixed points. 
\paragraph{Organization of the paper}
In the first section, we set up our notation, introduce basic notions and recall some results. In the 
second section, we compute bounds for the canonical Green's functions evaluated at two inequivalent cusps. 
{\small{\paragraph{Acknowledgements}
This article is part of the PhD thesis of the author, which was completed under the supervision 
of  J.~Kramer at Humboldt Universit\"at zu Berlin. The author would like to express his 
gratitude to J.~Kramer, J.~Jorgenson, and R.~S.~de~Jong for many interesting scientific discussions.}}  
\section{Background material}\label{section2}
Let $\Gamma \subset \mathrm{PSL}_{2}(\mathbb{R})$ be a Fuchsian subgroup of the first kind acting by 
fractional linear transformations on the upper half-plane $\mathbb{H}$. Let $X$ be the quotient space 
$\Gamma\backslash \mathbb{H}$, and let $g$ denote the genus of $X$. The quotient space $X$ admits the 
structure of a Riemann surface.  

\vspace{0.2cm}
Let $\mathcal{P}(\Gamma)$ and $H(\Gamma)$ denote the set of parabolic and hyperbolic elements of $\Gamma$, 
respectively. Let $\mathcal{P}$ be the finite set of cusps of $X$, respectively. Let $\overline{X}$ denote 
$\overline{X}=X\cup\mathcal{P}$. Locally, away from the cusps, we identity $\overline{X}$ with its universal 
cover $\mathbb{H}$, and hence, denote the points on $\overline{X}\backslash \mathcal{P}$ by the same letter as 
the points on $\mathbb{H}$.
\paragraph{Structure of $\overline{X}$ as a Riemann surface} 
The quotient space $\overline{X}$ admits the structure of a compact Riemann surface. We refer the reader to 
section 1.8 in $\cite{miyake}$, for the details regarding the structure of $\overline{X}$ as a compact 
Riemann surface. For the convenience of the reader, we recall the coordinate functions for the neighborhoods 
of cusps.

\vspace{0.2cm}
Let $p\in\mathcal{P}$ be a cusp and let $w\in U_{r}(p)$. Then $\vartheta_{p}(w)$ is given by
\begin{equation*}
 \vartheta_{p}(w)= e^{2\pi i \sigma_{p}^{-1}w},
\end{equation*}
where $\sigma_{p}$ is a scaling matrix of the cusp $p$ satisfying the following relations
\begin{align*}
\sigma_{p}i\infty = p \quad \mathrm{and} \quad \sigma_{p}^{-1}\Gamma_{p}\sigma_{p} = \langle\gamma_{\infty}\rangle,\quad
\mathrm{where}\,\,\, \gamma_{\infty}=\left(\begin{array}{ccc} 1 & 1\\ 0 & 1  \end{array}\right)
\quad&\mathrm{and}\quad\Gamma_{p}=\langle\gamma_{p}\rangle
\end{align*}
denotes the stabilizer of the cusp $p$ with generator $\gamma_{p}$.
\paragraph{Hyperbolic metric} 
We denote the (1,1)-form corresponding to the hyperbolic metric of $X$, which is compatible with the complex 
structure on $X$ and has constant negative curvature equal to minus one, by $\hyp(z)$. Locally, for 
$z\in X\backslash \mathcal{E}$, it is given by
\begin{equation*}
 \hyp(z)= \frac{i}{2}\cdot\frac{dz\wedge d\overline{z}}{{\Im(z)}^{2}}.
\end{equation*} 
Let $\vx(X)$ be the volume of $X$ with respect to the hyperbolic metric $\hyp$. It is given by the formula 
\begin{equation*}
\vx(X) = 2\pi\big(2g -2 + |\mathcal{P}|\big). 
\end{equation*}
The hyperbolic metric $\hyp(z)$ is singular at the cusps, and the rescaled hyperbolic metric
\begin{equation*}
 \shyp(z)= \frac{\hyp(z)}{ \vx(X)}
\end{equation*}
measures the volume of $X$ to be one. 

Locally, for $z$ $\in$ $X$, the hyperbolic Laplacian $\Delta_{\mathrm{hyp}}$ on  $X$ is given by
\begin{equation*}
\Delta_{\mathrm{hyp}} = -y^{2}\bigg(\frac{\partial^{2}}{\partial x^{2}} +
 \frac{\partial^{2}}{\partial y^{2}}\bigg) = -4y^{2}\bigg(\frac{\partial^{2}}{\partial z
\partial \overline{z}} \bigg). 
\end{equation*}
Recall that $d=\left(\partial + \overline{\partial} \right)$, $d^{c}=\dfrac{1}{4\pi i}\left( \partial - 
\overline{\partial}\right)$, and $dd^{c}= -\dfrac{\partial\overline{\partial}}{2\pi i}$. Furthermore, we 
have 
\begin{align}\label{ddcreln}
 d_{z}d_{z}^{c}=\del\hyp(z).
\end{align}
\paragraph{Canonical metric}
Let $S_{2}(\Gamma)$ denote the $\mathbb{C}$-vector space of cusp forms of 
weight 2 with respect to $\Gamma$ equipped with the Petersson inner product. Let 
$\lbrace f_{1},\ldots,f_{g}\rbrace $ denote an orthonormal basis of $S_{2}(\Gamma)$ with respect to the Petersson 
inner product. Then, the (1,1)-form $\can(z)$ corresponding to the 
canonical metric of $X$ is given by 
\begin{equation*}
 \can(z)=\frac{i}{2g} \sum_{j=1}^{g}\left|f_{j}(z)\right|^{2}dz\wedge d\overline{z}.
\end{equation*}
The canonical metric $\can(z)$ remains smooth at the cusps, and measures 
the volume of $X$ to be one. We denote the smooth (1,1)-form defined by $\can(z)$ on $\overline{X}$ by 
$\hatcan(z)$. 

\vspace{0.2cm}
For $z\in X$, we put,
\begin{align}\label{defndx}
d_{X}=\sup_{z\in X}\frac{\can(z)}{\shyp(z)}.
\end{align}
As the canonical metric $\can(z)$ remains smooth at the cusps and at the elliptic fixed points, 
and the hyperbolic metric is singular at these points, the quantity $d_{X}$ is well-defined. 
\paragraph{Canonical Green's function}
For $z, w \in X$, the canonical Green's function $\gcan(z,w)$ is defined as 
the solution of the differential equation (which is to be interpreted in terms of currents)
\begin{equation}\label{diffeqngcan}
d_{z}d^{c}_{z}\gcan(z,w)+ \delta_{w}(z)=\can(z),
\end{equation}
with the normalization condition
\begin{equation*}\label{normcondgcan}
 \int_{X}\gcan(z,w)\can(z)=0. 
\end{equation*}
From equation \eqref{diffeqngcan}, it follows that $\gcan(z,w)$ admits a $\log$-singularity at 
$z=w$, i.e., for $z, w\in X$, it satisfies 
\begin{equation}\label{gcanbounded}
\lim_{w\rightarrow z}\big(\gcan(z,w)+ \log |\vartheta_{z}(w)|^{2}\big)= O_{z}(1). 
\end{equation}
\paragraph{Parabolic Eisenstein Series}
For $z\in X$ and $s\in\mathbb{C}$ with $\Re(s)> 1$, the parabolic Eisenstein series 
$\mathcal{E}_{\mathrm{par},p}(z,s)$ corresponding to a cusp $p\in\mathcal{P}$ is defined by the 
series
\begin{equation*}
\mathcal{E}_{\mathrm{par},p}(z,s) = \sum_{\eta \in \Gamma_{p}\backslash \Gamma}
\Im(\sigma_{p}^{-1}\eta z)^{s}.
\end{equation*}
The series converges absolutely  and uniformly for $\Re(s) >1 $. It admits a meromorphic continuation to all 
$s\in\mathbb{C}$ with a simple pole at $s = 1$, and the Laurent expansion at $s=1$ is of the form 
\begin{equation}\label{laurenteisenpar}
\mathcal{E}_{\mathrm{par},p}(z,s) = \frac{1}{(s-1)\vx(X)} + \kappa_{p}(z) + O_{z}(s-1),
\end{equation}
where $\kappa_{p}(z)$ the constant term of $\mathcal{E}_{\mathrm{par},p}(z,s)$ at $s=1$ is called 
Kronecker's limit function (see Chapter 6 of \cite{hi}).
\vspace{0.2cm}
For $z\in X$, and $p,q \in \mathcal{P}$, the Kronecker's limit function 
$\kappa_{p}(\sigma_{q}z)$ satisfies the following equation (see Theorem 1.1 of \cite{jsu} for the proof)
\begin{align}\label{fourierkappaeqn}
\kappa_{p}(\sigma_{q}z)= \sum_{n < 0} k_{p,q}(n)e^{2 
\pi in\overline{z}}+ \delta_{p,q}\Im(z)+
k_{p,q}(0)- \frac{\log\big(\Im(z)\big)}{\vx(X)} + 
\sum_{n > 0}k_{p,q}(n)e^{2\pi i nz},   
\end{align}
with Fourier coefficients $k_{p,q}(n)$ $\in$ $\mathbb{C}$.
Let $p,q \in \mathcal{P}$, then for $z\in X$ and $s\in\mathbb{C}$ with $\Re(s)> 1$, 
the parabolic Eisenstein series $\mathcal{E}_{\mathrm{par,}p}(\sigma_{q}z,s)$ 
associated to $p\in \mathcal{P}$, admits a Fourier expansion of the form
\begin{align}\label{foureireisenstein}
\mathcal{E}_{\mathrm{par},p}(\sigma_{q}z,s)= \delta_{p,q}y^{s} + \alpha_{p,q}(s)y^{1-s} + 
\sum_{n\not=0}\alpha_{p,q}(n,s)W_{s}(nz),
\end{align}
where $\alpha_{p,q}(s)$, $\alpha_{p,q}(n,s)$, and $W_{s}(nz)$ the Whittaker function are given by equations 
(3.21), (3.22), and (1.37), respectively in \cite{hi}. 
\paragraph{Heat Kernels}
For $t \in \mathbb{R}_{> 0}$  and $z, w \in \mathbb{H}$, the hyperbolic heat kernel $K_{\mathbb{H}}(t;z,w)$ on $\mathbb{R}_{> 0}\times
\mathbb{H} \times\mathbb{H}$ is given by the formula
\begin{equation*}\label{defnkh}
K_{\mathbb{H}}(t;z,w)= \frac{\sqrt{2}e^{- t\slash 4}}{(4\pi t)^{3\slash 2}}
\int_{d_{\mathbb{H}}(z,w)}^{\infty}\frac{re^{-r^{2}\slash 4t}}{\sqrt{\cosh(r)-\cosh (d_{\mathbb{H}}(z,w))}}dr,
\end{equation*}
where $d_{\mathbb{H}}(z,w)$ is the hyperbolic distance between $z$ and $w$.

For $t \in  \mathbb{R}_{> 0}$ and $z, w \in X$, the hyperbolic heat kernel $\khyp(t;z,w)$ on $\mathbb{R}_{> 0}\times X\times X$ is defined as 
\begin{equation*}\label{defnkhyp}
\khyp(t;z,w)=\sum_{\gamma\in\Gamma}K_{\mathbb{H}}(t;z,\gamma w).
\end{equation*}
For $z,w\in X,$ the hyperbolic heat kernel $\khyp(t;z,w)$ satisfies the differential equation
\begin{align}\label{diffeqnkhyp} 
\bigg(\Delta_{\text{hyp},z} + \frac{\partial}{\partial t}\bigg)\khyp(t;z,w) &=0,
\end{align}
Furthermore for a fixed $w\in X$ and any smooth function $f$ on $X$, the hyperbolic heat kernel $\khyp(t;z,w)$ 
satisfies the equation
\begin{align}\label{normcondkhyp}  
\lim_{t\rightarrow 0}\int_{X}\khyp(t;z,w)f(z)\hyp(z) &= f(w).
\end{align}
To simplify notation, we write $\khyp(t;z)$ instead of $\khyp(t;z,z)$, when $z=w$.

\vspace{0.2cm}
For $t\in \mathbb{R}_{\geq 0}$ and $z\in X$, put 
\begin{align*}
\hkhyp(t;z)=\sum_{\gamma\in\mathcal{P}(\Gamma)}\kh(t;z,\gamma z).
\end{align*}
The convergence of the above series follows from the convergence of the hyperbolic heat kernel 
$\khyp(t;z)$ and the fact that $\kh(t;z,\gamma z)$ is positive for all $t\in\mathbb{R}_{\geq 0}$, $z\in 
\mathbb{H}$, and $\gamma\in\Gamma$. 
\paragraph{Selberg constant}
The hyperbolic length of the closed geodesic determined by a primitive non-conjugate hyperbolic element 
$\gamma\in \mathcal{H}(\Gamma)$ on $X$ is given by 
\begin{align*}
\ell_{\gamma}=\inf\lbrace{d_{\mathbb{H}}(z,\gamma z)|\,z\in\mathbb{H}\rbrace}.
\end{align*}

\vspace{0.2cm}
For $s\in\mathbb{C}$ with $\Re(s)>1$, the Selberg zeta function associated to $X$ is defined as
\begin{equation*}
Z_{X}(s)= \prod_{\gamma\in \mathcal{H}(\Gamma)}Z_{\gamma}(s), \,\,\,\,\,\,\mathrm{where}\,\,\,Z_{\gamma}(s)= 
\prod_{n=0}^{\infty}\big(1-e^{(s+n)\ell_{\gamma}}\big).
\end{equation*}
The Selberg zeta function $Z_{X}(s)$ admits a meromorphic continuation to all $s\in\mathbb{C}$, with 
zeros and poles characterized by the spectral theory 
of the hyperbolic Laplacian. Furthermore, $Z_{X}(s)$ has a simple zero at 
$s=1$, and the following constant is well-defined
\begin{equation}\label{defnselbergzeta}
c_{X}= \lim_{s\rightarrow 1}\bigg(\frac{Z^{'}_{X}(s)}{Z_{X}(s)}-
\frac{1}{s-1}\bigg).
\end{equation}
For $t\in\mathbb{R}_{\geq 0}$, the hyperbolic heat trace is given by the integral
\begin{align*}
H\mathrm{Tr}\khyp(t)=\int_{X}\hkhyp(t;z)\hyp(z). 
\end{align*}
The convergence of the integral follows from the celebrated Selberg trace formula. Furthermore, from Lemma 4.2 in \cite{jk1}, 
we have the following relation    
\begin{align}\label{selbergconstant}
\int_{0}^{\infty} \big(H\mathrm{Tr}\khyp(t)-1\big)dt=c_{X}-1.
\end{align}
\paragraph{Automorphic Green's function}
For $z, w \in \mathbb{H}$ with $z \not= w$, and  $s$ $\in$ $\mathbb{C}$ with $\Re(s)> 0$, the free-space Green's function $g_{\mathbb{H},s}(z,w)$ is defined as
\begin{equation*} 
g_{\mathbb{H},s}(z,w) = g_{\mathbb{H},s}(u(z,w))= \dfrac{\Gamma(s)^{2}}{\Gamma(2s)}u^{-s}
F(s,s;2s,-1\slash u),
\end{equation*} 
where $u=u(z,w)=|z-w|^{2}\slash( 4\Im(z)\Im(w))$ and $F(s,s;2s,-1\slash u)$ is the 
hypergeometric function. 

There is a sign error in the formula defining the free-space Green's function given by equation (1.46) in \cite{hi}, i.e., 
the last argument $-1\slash u$ in the hypergeometric function has been incorrectly stated as $1\slash u$, which we have 
corrected in our definition. We have also normalized the free-space Green's function defined in \cite{hi} by multiplying it by 
$4\pi.$

\vspace{0.2cm}
For $z, w \in X$ with $z\not = w$, and $s\in\mathbb{C}$ with $\Re(s) > 1$, the automorphic Green's function 
$g_{\mathrm{hyp},s}(z,w)$ is defined as
\begin{equation*}
g_{\mathrm{hyp},s}(z,w) = \sum_{\gamma\in\Gamma}g_{\mathbb{H},s}(z,\gamma w).
\end{equation*}
The series converges absolutely uniformly for $z\not = w$ and $\Re(s) > 1$ (see Chapter 5 in \cite{hi}). 

For $z, w \in X$ with $z \not = w$, and $s\in\mathbb{C}$ with $\Re(s) > 1$, the automorphic Green's function satisfies the 
following properties (see Chapters 5 and 6 in \cite{hi}):

(1) For $\Re(s(s-1)) > 1$, we have
\begin{equation*}
g_{\mathrm{hyp},s}(z,w) = 4\pi\int_{0}^{\infty}\khyp(t;z,w)e^{-s(s-1)t}dt.
\end{equation*}
(2) It admits a logarithmic singularity along the diagonal, i.e.,  
\begin{equation*}
\lim_{w\rightarrow z}\big(g_{\mathrm{hyp},s}(z,w) + \log{|\vartheta_{z}(w)|^{2}}\big)= O_{s,z}(1).
\end{equation*}
(3) The automorphic Green's function $g_{\mathrm{hyp},s}(z,w)$ admits a meromorphic continuation to all $s\in\mathbb{C}$ with a 
simple pole at $s=1$ with residue $4\pi\slash\vx(X)$, and the Laurent expansion at $s=1$ is of the form
\begin{equation*}
g_{\mathrm{hyp},s}(z,w)= \frac{4\pi}{s(s-1)\vx(X)} + g^{(1)}_{\mathrm{hyp}}(z,w) + O_{z,w}(s-1),
\end{equation*} 
where $g_{\mathrm{hyp}}^{(1)}(z,w)$ is the constant term of $g_{\mathrm{hyp},s}(z,w)$ at $s=1$. 

\vspace{0.15cm}
(4) Let $p,q\in \mathcal{P}$ be two cusps. Put
 \begin{align*}
C_{p,q}  = \min \bigg{\lbrace} c > 0\,\bigg{|} \bigg(\begin{array}{ccc} a &b\\
 c & d  \end{array}\bigg) \in \sigma_{p}^{-1}\Gamma \sigma_{q} \bigg{\rbrace},
\end{align*} 
and $C_{p,p}=C_{p}$. Then, for $z,w\in X$ with $\Im(w) > \Im(z)$ and  $\Im(w)\Im(z) > C_{p,q}^{-2}$, and $s\in\mathbb{C}$ with $\Re(s) > 1$, 
the automorphic Green's function admits the Fourier expansion
\begin{align} 
g_{\mathrm{hyp},s}(\sigma_{p}z,\sigma_{q}w)= \frac{4\pi\Im(w)^{1-s}}{2s-1}
\mathcal{E}_{\mathrm{par},q}(\sigma_{p}z,s) -\delta_{p,q} \log\big|1- e^{2\pi i (w-z)}\big|^{2} 
+O\big(e^{-2\pi (\Im( w)-\Im (z))}\big). \label{fourierautghyp}
\end{align}
This equation has been proved as Lemma 5.4 in \cite{hi}, and one of the terms was wrongly estimated in the proof of the lemma. 
We have corrected this error, and stated the corrected equation. 
\paragraph{The space $C_{\ell,\ell\ell}(X)$} 
Let $C_{\ell,\ell\ell}(X)$ denote the set of complex-valued functions $f:X\rightarrow \mathbb{P}^{1}(
\mathbb{C})$, which admit the following type of singularities at finitely many points $\mathrm{Sing}(f)\subset X$, 
and are smooth away from $\mathrm{Sing}(f)$: 

\vspace{0.2cm}
(1) If $s\in\mathrm{Sing}(f)$, then as $z$ approaches $s$, the function $f$ satisfies  
\begin{align}\label{fsingular}
f(z)= c_{f,s}\log|\vartheta_{s}(z)| + O_{z}(1),
\end{align}
for some $c_{f,s}\in \mathbb{C}$.

\vspace{0.2cm}
(2) As $z$ approaches a cusp $p\in \mathcal{P}$, the function $f$ satisfies
\begin{align}\label{fcusp}
f(z)=c_{f,p}\log\big(-\log|\vartheta_{p}(z)|\big) + O_{z}(1),
\end{align}
for some $c_{f,p}\in \mathbb{C}$. 
\paragraph{Hyperbolic Green's function}
For $z, w \in X$ and $z\not = w$, the hyperbolic Green's function is defined as 
\begin{equation*}
\ghyp(z,w) = 4\pi\int_{0}^{\infty}\bigg(\khyp(t;z,w)-\frac{1}{\vx(X)}\bigg)dt.
\end{equation*}
For $z, w \in X$ with $z \not = w$, the hyperbolic Green's function satisfies the 
following properties:

(1) For $z, w \in X$, we have 
\begin{equation}\label{ghypbounded}
\lim_{w\rightarrow z}\big( \ghyp(z,w) + \log{|\vartheta_{z}(w)|^{2}}\big)= O_{z}(1).
\end{equation}
(2) For $z, w \in X$, the hyperbolic Green's function satisfies 
the differential equation (which is to be interpreted in terms of currents)
\begin{align}
d_{z}d_{z}^{c}\ghyp(z,w) +\delta_{w}(z)& = \shyp(z), \label{diffeqnghyp} \\
\intertext{ with the normalization condition} 
\int_{X}\ghyp(z,w)\hyp(z) & = 0. \label{normcondghyp}
\end{align}
(3) For $z,w\in X$ and $z\not=w$, we have
\begin{equation}\label{laurentghyp}
 \ghyp(z,w)= g^{(1)}_{\mathrm{hyp}}(z,w)= \lim_{s\rightarrow 1}\bigg(g_{\text{hyp},s}(z,w) - 
\frac{4\pi}{s(s-1)\vx(X)}\bigg).
\end{equation}
The above properties follow from the properties of the heat kernel $\khyp(t;z,w)$ (equations \eqref{diffeqnkhyp} and \eqref{normcondkhyp}) or 
from that of the automorphic Green's function $g_{\mathrm{hyp},s}(z,w)$. 

\vspace{0.15cm}
(4) From Proposition 2.4.1 in \cite{anilthesis}, (or from Proposition 2.1 in \cite{anilpaper1}) 
for a fixed $ w\in X$, and for $z\in X$ with $\Im(\sigma_{p}^{-1}z)>\Im(\sigma_{p}^{-1} w)$, and 
$\Im(\sigma_{p}^{-1}z)\Im(\sigma_{p}^{-1}w) >C_{p}^{-2}$, we have
\begin{align}
& \ghyp(z,w) = 4\pi\kappa_{p}(w) - 
\frac{4\pi}{\vx(X)}-\frac{4\pi\log\big(\Im(\sigma_{p}^{-1}z)\big)}{\vx(X)}-\notag \\ &
 \log\big{|}1-e^{2\pi i(\sigma_{p}^{-1}z - \sigma_{p}^{-1}w)}\big{|}^{2}+
O\big(e^{-2\pi (\Im(\sigma_{p}^{-1}z)-\Im(\sigma_{p}^{-1}w))}\big),\label{ghypcusp}
\end{align}
i.e., for a fixed $w\in X$, as $z\in X$ approaches a cusp $p\in\mathcal{P}$, we have
\begin{align}\label{ghyploglog}
\ghyp(z,w) & = -\frac{4\pi\log\big(\Im(\sigma_{p}^{-1}z)\big)}{\vx(X)}+
O_{z,w}(1) = -\frac{4\pi\log\big(-\log|\vartheta_{p}(z)|\big)}{\vx(X)}+ O_{z,w}(1).
\end{align}
(5) For any $f\in C_{\ell,\ell\ell}(X)$ and for any fixed $w\in X\backslash \mathrm{Sing}(f)$, 
from  Corollary 3.1.8  in \cite{anilthesis} (or from Corollary 2.5 in \cite{anilpaper1}), we have the 
equality of integrals 
\begin{align}\label{ghypcurrent}
\int_{X}\ghyp(z,w)d_{z}d_{z}^{c}f(z) + f(w)+
\sum_{s\in \mathrm{Sing}(f)}
\frac{c_{f,s}}{2}\ghyp(s,w)= \int_{X}f(z)\shyp(z).
\end{align}
\paragraph{Certain Convergence results} 
 For $z\in \mathbb{H}$, put
\begin{align*}
 P(z)= \sum_{\gamma\in\mathcal{P}(\Gamma)}\gh(z,\gamma z).
\end{align*}
The above series is invariant under the action $\Gamma$, and hence, defines a function on $X$. Furthermore, 
from Proposition 4.2.4 in \cite{anilthesis} (or from 2.2 in \cite{anilpaper2}), the above series 
converges for all $z\in X$, and satisfies the following equation
\begin{align} 
P(z)=\sum_{p\in \mathcal{P}}\sum_{\eta\in\Gamma_{p}\backslash\Gamma}P_{\mathrm{gen},p}(\eta z),
\label{decomposition}
\end{align}
where $P_{\mathrm{gen},p}(z)=\displaystyle\sum_{ n\not = 0}g_{\mathbb{H}}(z,\gamma_{p}^{n}z)$. 

\vspace{0.2cm}
Furthermore, from the absolute and uniform convergence of $P(z)$, and from that of the following series 
from Lemma 5.2 in \cite{jk}
\begin{align*}
\sum_{\gamma\in\mathcal{P}(\Gamma)}\del g_{\mathbb{H}}(z,\gamma z),
\end{align*}
we get
\begin{align}
&\sum_{\gamma\in\mathcal{P}(\Gamma)}\del g_{\mathbb{H}}(z,\gamma z)  =\del P(z)= 
\sum_{p\in \mathcal{P}}\sum_{\eta\in \Gamma_{p}\backslash \Gamma}\del P_{\mathrm{gen},p}(
\eta z),\label{delP1}\\& \del P_{\mathrm{gen},p}(z)=\sum_{n\not = 0}\del g_{\mathbb{H}}(\sigma_{p}^{-1} z,\gamma_{\infty}^{n}
\sigma_{p}^{-1}z)= 2\bigg(\frac{2\pi \Im(\sigma_{p}^{-1}z)}{\sinh(2\pi
\Im(\sigma_{p}^{-1}z))}\bigg)^{2}-2.\label{delP2}
\end{align}
For $z\in X$, put
\begin{equation}\label{H(z)defneqn}
H(z)=4\pi\int_{0}^{\infty}\bigg(\hkhyp(t;z)-\frac{1}{\vx(X)}\bigg)dt.  
\end{equation}
The function $H(z)$ is invariant under the action of $\Gamma$, and hence, defines a function on $X$. Furthermore, 
from Proposition 4.3.2 (or from Proposition 2.9), it follows that $H(z)$ is well-defined on $X$, and 
for $z,w\in X$, we have 
\begin{align*}
H(z)=  \lim_{w\rightarrow z}\big{(}\ghyp(z,w)-g_{\mathbb{H}}(z,w)\big{)}-P(z).
\end{align*}
From the above equation, for $z\in X$, we find
\begin{align*}
\del P(z)+ \del H(z)= \del\lim_{w\rightarrow z}\big(g_{\mathrm{hyp}}(z,w)-g_{\mathbb{H}}(z,w)\big).
\end{align*}
For $z\in X$, since the integral
\begin{align*}
4\pi\int_{0}^{\infty}\bigg(\khyp(t;z,z)-K_{\mathbb{H}}(t;0)-\frac{1}{\vx(X)}\bigg)dt,
\end{align*}
as well as the integral of the derivatives of the integrand are absolutely 
convergent, we can take the Laplace operator $\del$ inside the integral. So for $z\in X$, we find
\begin{align}\label{delkdecomposition}
\del P(z)+\del H(z)= 4\pi\int_{0}^{\infty}\del K_{\mathrm{hyp}}(t;z)dt.
\end{align}
From Lemma 5.2 and Proposition 7.3 in \cite{K}, for $z\in X$, we have the following relation 
\begin{align*}
4\pi\int_{0}^{\infty}\del \khyp(t;z)dt=\sum_{\gamma\in\Gamma\backslash \lbrace\id\rbrace}\gh(z,\gamma z)
\end{align*}
and the right-hand side of above equation remains bounded at the cusps. So we deduce that the left-hand side 
also remains bounded at the cusps. 

\vspace{0.2cm}
From  Proposition 4.3.3 in \cite{anilthesis} (or from Proposition 2.10 in \cite{anilpaper2}), 
for $z\in X$ approaching a cusp $p\in\mathcal{P}$, we have
\begin{align}\label{H(z)cusp}
H(z)&\,= -\frac{8\pi}{\vx(X)}\log\big(\Im(\sigma_{p}^{-1}z)\big)-\frac{4\pi}{\vx(X)}+
4\pi k_{p,p}(0)+O\big(\Im(\sigma_{p}^{-1}z)^{-1}\big)\notag\\&\,=-\frac{8\pi}{\vx(X)}\log
\big(-\log|\vartheta_{p}(z)|\big)+O_{z}(1). 
\end{align}
Hence, we can conclude that the function $H(z)\in C_{\ell,\ell\ell}(X)$ with $\mathrm{Sing}(f)=\emptyset$. 
Lastly, from equation \eqref{selbergconstant}, we have
\begin{align}\label{H(z)integral}
\int_{X}H(z)\hyp(z)=4\pi (c_{X}-1).
\end{align}
\paragraph{An auxiliary identity}
For notational brevity, put
\begin{align*}
C_{\mathrm{hyp}}= \int_{X}\int_{X}\ghyp(\zeta,\xi)\bigg(\int_{0}^{\infty}\del 
\khyp(t;\zeta)dt\bigg)\bigg(\int_{0}^{\infty}\del \khyp(t;\xi)dt\bigg)\hyp(\xi)\hyp(\zeta).
\end{align*}
From Proposition 2.6.4 in \cite{anilthesis} (or from Proposition 2.8 in \cite{anilpaper1}, for 
$z,w\in X$, we have
\begin{equation}\label{phi}
\ghyp(z,w)-\gcan(z,w)= \phi(z) + \phi(w),
\end{equation}
where from Corollary 3.2.7 in \cite{anilthesis} (or from  Remark 2.16 in \cite{anilpaper1}), 
the function $\phi(z)$ is given by the formula 
\begin{align*}
&\phi(z)= \frac{1}{2g}\int_{X}\ghyp(z,\zeta)
\left(\int_{0}^{\infty}\del \khyp(t;\zeta)dt\right)\hyp(\zeta)-\frac{C_{\mathrm{hyp}}}{8g^{2}}.
\end{align*} 
As $H(z)\in C_{\ell,\ell\ell}(X)$, using the relations \eqref{ddcreln} and \eqref{H(z)integral}, and 
combining equation \eqref{delkdecomposition} with \eqref{ghypcurrent}, we derive
\begin{align}\label{phi(z)formula}
\phi(z)&\,=\frac{H(z)}{2g}+\frac{1}{8\pi g}\int_{X}g_{\mathrm{hyp}}(z,\zeta)\del P(\zeta)
\hyp(\zeta) - \frac{C_{\mathrm{hyp}}}{8g^{2}}-\int_{X}H(z)\shyp(z)\notag\\&\,=
\frac{H(z)}{2g}+\frac{1}{8\pi g}\int_{X}g_{\mathrm{hyp}}(z,\zeta)\del P(\zeta)
\hyp(\zeta) - \frac{C_{\mathrm{hyp}}}{8g^{2}}-\frac{2\pi(c_{X}-1)}{g\vx(X)}. 
\end{align}
For more details regarding the proof for the above computation, we refer the reader to Theorem 4.3.8 in 
\cite{anilthesis} (or Corollary 2.12 in \cite{anilpaper2}).
\paragraph{Key identity}
For $z \in X$, we have the relation of differential forms
\begin{align*}
&g\can(z) =\bigg(\frac{1}{4\pi}+ \frac{1}{\vx(X)}\bigg)\hyp(z)+ \frac{1}{2}\bigg(\int_{0}^{\infty}\del
\khyp(t;z) dt\bigg)\hyp(z).
\end{align*}
This relation has been established as Theorem 3.4 in \cite{jk}, when $X$ is 
compact, which easily extends to our case. Furthermore, from Corollary 3.2.5 in \cite{anilthesis} 
(or from Corollary 2.15 in \cite{anilpaper1}), for any $f\in C_{\ell,\ell\ell}(X)$, we have
\begin{align}\label{keyidentity}
& g\int_{X}f(z)\can(z) = \notag\\&\bigg(\frac{1}{4\pi}+\frac{1}{\vx(X)} \bigg)
\int_{X}f(z)\hyp(z) + \frac{1}{2}\int_{X}f(z)\bigg(\int_{0}^{\infty}
\del \khyp(t;z)dt \bigg)\hyp(z).  
\end{align}
\section{Bounds for canonical Green's functions at cusps}
Let $p,q \in \mathcal{P}$ be two cusps with $p\not =q $. Then, from equation \eqref{phi}, 
we find
\begin{align}
g_{\mathrm{can}}(p,q)= \lim_{z\rightarrow p}\lim_{w\rightarrow q}
\big(g_{\mathrm{hyp}}(z,w)-\phi(z)-\phi(w)\big).
\end{align}
From equations \eqref{ghypcusp} and \eqref{H(z)cusp}, we know the asymptotics of the functions 
$\ghyp(z,w)$ and $H(z)$ at the cusps, respectively. So if we can compute the asymptotics of the integral
\begin{align*}
\int_{X}g_{\mathrm{hyp}}(z,\zeta)\del P(\zeta)\hyp(\zeta)
\end{align*}
at the cusps, we will be able to compute an upper bound for the canonical 
Green's function when evaluated at two different cusps. 

\vspace{0.2cm}\noindent 
For the remaining part of the thesis, for $p\in \mathcal{P}$ a cusp and 
$z\in\mathbb{H}$, we denote $\Im(\sigma_{p}^{-1}z)$ by $y_{p}.$ 

\vspace{0.2cm}\noindent 
In the following two lemmas, we compute the zeroth Fourier coefficient of the automorphic Green's function 
and the hyperbolic Green's function.  
\begin{lem}\label{lem3.1}
Let $p,q\in\mathcal{P}$ be two cusps. Then, for $z$ and $w=u+iv\in X$ with 
$y_{p}> v $ and $v y_{p}> 1$, and $s\in\mathbb{C}$ with $\Re(s)> 1$, we have 
\begin{align}
 \int_{0}^{1}g_{\mathrm{hyp},s}(z,\sigma_{q}w)du=\frac{4\pi v^{1-s}}{2s-1}\mathcal{E}_{\mathrm{par},q}(z,s)+
\frac{ 4\pi\delta_{p,q}}{2s-1}\big(v^{s}y_{p}^{1-s}-v^{1-s}y_{p}^{s}\big).
\label{zerocoeffecient1automorphic}
\end{align}
Furthermore, for $v > y_{p} $ and $v y_{p}> 1$, and $s\in\mathbb{C}$ with $\Re(s)> 1$, we have
\begin{align}\label{zerocoeffecient2automorphic}
\int_{0}^{1}g_{\mathrm{hyp},s}(z,\sigma_{q}w)du=\frac{4\pi v^{1-s}}{2s-1}\mathcal{E}_{\mathrm{par},q}(z,s).
\end{align} 
\begin{proof}
For $z$ and $w=u+iv\in X$ with $y_{p}> v $ and $v y_{p}> 1$, and $s\in\mathbb{C}$ with $\Re(s)> 1$, 
combining Lemmas 5.1 and 5.2 of \cite{hi}, we have 
\begin{align*}
\int_{0}^{1}g_{\mathrm{hyp},s}(z,\sigma_{q}w)du= \frac{4\pi y_{p}^{1-s}}{2s-1}\left(\delta_{p,q}
v^{s} + \alpha_{p,q}(s)v^{1-s}\right)+\frac{4\pi v^{1-s}}{2s-1}\sum_{n\not = 0}
\alpha_{p,q}(n,s)W_{s}(n\sigma_{p}^{-1}z).
\end{align*}
The expression on the right-hand side of the above equation can be rewritten as
\begin{align}
&\frac{4\pi v^{1-s}}{2s-1}\bigg(\delta_{p,q}y_{p}^{s}+\alpha_{p,q}(s)y_{p}^{1-s}+
\sum_{n\not = 0}\alpha_{p,q}(n,s)W_{s}(n\sigma_{p}^{-1}z)\bigg)+
\frac{ 4\pi\delta_{p,q}}{2s-1}\big(v^{s}y_{p}^{1-s}-v^{1-s}y_{p}^{s}\big).\label{autgreen01}
\end{align}
For $s\in\mathbb{C}$ and $\Re(s)> 1$, from the Fourier expansion of the parabolic Eisenstein series 
$\mathcal{E}_{\mathrm{par},q}(z,s)$ described in equation \eqref{foureireisenstein}, we get 
\begin{align}
\frac{4\pi v^{1-s}}{2s-1}\bigg(\delta_{p,q}y_{p}^{s}+\alpha_{p,q}(s)y_{p}^{1-s}+\sum_{n\not = 0}\alpha_{p,q}
(n,s)W_{s}(n\sigma_{p}^{-1}z)\bigg)=\frac{4\pi v^{1-s}}{2s-1}\mathcal{E}_{\mathrm{par},q}(z,s).\label{autgreen02}
\end{align}
Combining equations (\ref{autgreen01}) and (\ref{autgreen02}) proves equation (\ref{zerocoeffecient1automorphic}). 

\vspace{0.2cm}\noindent 
For $v > y_{p} $ and $v y_{p}> 1$, and $s\in\mathbb{C}$ with $\Re(s)> 1$, combining Lemmas 5.1 and 5.2 of \cite{hi}, we have 
\begin{align*}
&\int_{0}^{1}g_{\mathrm{hyp},s}(z,\sigma_{q}w)du=\frac{4\pi v^{1-s}}{2s-1}\bigg(\delta_{p,q}y_{p}^{s}+\alpha_{p,q}(s)y_{p}^{1-s}
+\sum_{n\not = 0}\alpha_{p,q}(n,s)W_{s}(n\sigma_{p}^{-1}z)\bigg).
\end{align*}
From equation (\ref{autgreen02}), we derive that
\begin{align*}
 \int_{0}^{1}g_{\mathrm{hyp},s}(z,\sigma_{q}w)du=\frac{4\pi v^{1-s}}{2s-1}\mathcal{E}_{\mathrm{par},q}(z,s),
\end{align*}
which proves equation (\ref{zerocoeffecient2automorphic}), and completes the proof of the lemma.
\end{proof}
\end{lem}
\begin{lem}\label{lem3.2}
Let $p,q\in\mathcal{P}$ be two cusps. Then, for $z$ and $w=u+iv\in X$ with 
$y_{p}> v $ and $v y_{p}> 1$, we have 
\begin{align}
\int_{0}^{1}g_{\mathrm{hyp}}(z,\sigma_{q}w)du= 4\pi\kappa_{q}(z)-\frac{4\pi}{\vx(X)}-\frac{4\pi\log v}{\vx(X)}+
4\pi\delta_{p,q}(v-y_{p}).\label{zerocoeffecient1}
\end{align}
Furthermore, for $v > y_{p} $ and $v y_{p}> 1$, we have
\begin{align}\label{zerocoeffecient2}
\int_{0}^{1}g_{\mathrm{hyp}}(z,\sigma_{q}w)du= 4\pi\kappa_{q}(z)-\frac{4\pi}{\vx(X)}-
\frac{4\pi\log v}{\vx(X)}.
\end{align}
\begin{proof}
Observe that
\begin{align}
&\int_{0}^{1}g_{\mathrm{hyp}}(z,\sigma_{q}w)du=\int_{0}^{1}\lim_{s\rightarrow 1}\bigg(g_{\mathrm{hyp},s}
(z,\sigma_{q}w)-\frac{4\pi}{s(s-1)\vx(X)}\bigg)du=\notag\\&\lim_{s\rightarrow 1}\bigg(\int_{0}^{1}
g_{\mathrm{hyp},s}(z,\sigma_{q}w)du-\frac{4\pi}{(s-1)\vx(X)}\bigg)+\frac{4\pi}{\vx(X)}.
\label{zerofouriercoeffecienteqn4}
\end{align}
For $z$ and $w=u+iv\in X$ with $y_{p}> v $ and $v y_{p}> 1$, combining equations (\ref{zerocoeffecient1automorphic}) 
and (\ref{zerofouriercoeffecienteqn4}), we find that the right-hand side of the above equation decomposes into 
the following expression
\begin{align*}
\lim_{s\rightarrow 1}\bigg(\frac{4\pi v^{1-s}}{2s-1}
\mathcal{E}_{\mathrm{par},q}(z,s)-\frac{4\pi}{(s-1)\vx(X)}\bigg) +4\pi\delta_{p,q}(v-y_{p})+
\frac{4\pi}{\vx(X)}.
\end{align*}
To evaluate the above limit, we compute the Laurent expansions of $\mathcal{E}_{\mathrm{par},p}(w,s)$, 
$\Im(\sigma_{p}^{-1}z)^{1-s}$, and $(2s-1)^{-1}$ at $s=1$. The Laurent expansions of 
$\Im{(\sigma_{p}^{-1}z)^{1-s}}$ and $(2s-1)^{-1}$ at $s=1$ are easy to compute, and are of the form
\begin{align*}
&\Im{(\sigma_{p}^{-1}z)}^{1-s}=1 - (s-1)\log\big(\Im{(\sigma_{p}^{-1}z)}\big)
+ O\big((s-1)^{2}\big); \,\,\,\,\frac{1}{2s-1}  = 1- 2(s-1) + O\big((s-1)^{2}\big).
\end{align*}
Combining the above two equations with equation \eqref{laurenteisenpar}, we find 
\begin{align}\label{zerolim}
&4\pi\lim_{s\rightarrow 1}\bigg(\frac{\Im(\sigma_{p}^{-1}z)^{1-s}}{2s-1}
\mathcal{E}_{\mathrm{par},p}(w,s)-\frac{1}{(s-1)\vx(X)}\bigg)=\notag\\&
 4\pi \kappa_{p}(w) -\frac{8\pi}{\vx(X)}- \frac{4\pi\log\big(\Im(\sigma_{p}^{-1}z)\big)}{\vx(X)},
\end{align}
Combining the above computation with equation \eqref{zerofouriercoeffecienteqn4}, we arrive at
\begin{align*}
&\int_{0}^{1}g_{\mathrm{hyp}}(z,\sigma_{q}w)du= 4\pi\kappa_{q}(z)-
\frac{4\pi}{\vx(X)}-\frac{4\pi\log v}{\vx(X)}+4\pi\delta_{p,q}(v-y_{p}),
\end{align*}
which proves equation (\ref{zerocoeffecient1}).

\vspace{0.2cm}\noindent We now prove equation (\ref{zerocoeffecient2}). For $v > y_{p} $ and $v y_{p}> 1$, 
combining equations (\ref{zerocoeffecient2automorphic}) and (\ref{zerofouriercoeffecienteqn4}), we find
\begin{align}
&\int_{0}^{1}g_{\mathrm{hyp}}(z,\sigma_{q}w)du=\lim_{s\rightarrow 1}\bigg(\frac{4\pi v^{1-s}}{2s-1}
\mathcal{E}_{\mathrm{par},q}(z,s)-\frac{4\pi}{(s-1)\vx(X)}\bigg) +\frac{4\pi}{\vx(X)}.
\label{zerofouriercoeffecienteqn5}
\end{align}
Combining equations (\ref{zerofouriercoeffecienteqn5}) and (\ref{zerolim}), we find 
\begin{align*}
&\int_{0}^{1}g_{\mathrm{hyp}}(z,\sigma_{q}w)du=4\pi\kappa_{q}(z)-\frac{4\pi}{\vx(X)}-\frac{4\pi\log v}{\vx(X)},
\end{align*}
which proves equation (\ref{zerocoeffecient2}), and hence, completes the proof of the lemma.
\end{proof}
\end{lem}
\begin{prop}\label{propghypdelpint}
Let $p\in\mathcal{P}$ be a cusp. For $z$, $w=u+iv\in X$ with $y_{p}>1$, we have the formal decomposition
\begin{align}
&\int_{X}g_{\mathrm{hyp}}(z,w)\del P(w)\hyp(w)=\sum_{q\in \mathcal{P}}\int_{0}^{1\slash y_{p}}
\int_{0}^{1}g_{\mathrm{hyp}}(z,\sigma_{q}w)\del P_{\mathrm{gen},q}(\sigma_{q}w)\frac{dudv}{v^{2}}+\notag\\
&\sum_{q\in \mathcal{P}}\int_{1\slash y_{p}}^{\infty}\bigg(4\pi\kappa_{q}(z)-\frac{4\pi}{\vx(X)}-
\frac{4\pi\log v}{\vx(X)}\bigg)\del P_{\mathrm{gen},q}(\sigma_{q}w) \frac{dv}{v^{2}}+\notag\\&
4\pi\int_{1\slash y_{p}}^{y_{p}} (v-y_{p})\del P_{\mathrm{gen},p}
(\sigma_{p}w)\frac{dv}{v^{2}}.
\label{propghypdelpinteqn}
\end{align}
\begin{proof}
As the series $\del P(w)$ is absolutely and uniformly convergent, we have
\begin{align}
 &\int_{X}g_{\mathrm{hyp}}(z,w)\del P(w)\hyp(w)=
\sum_{q\in \mathcal{P}}\sum_{\eta\in\Gamma_{q}\backslash \Gamma}
\int_{X}g_{\mathrm{hyp}}(z,w)\del P_{\mathrm{gen},q}(\eta w)\hyp(w),\label{gcan(pq)eqn1}
\end{align}
After making the substitution $w\mapsto \eta^{-1}\sigma_{q}w$, from the $\Gamma$-invariance of $g_{\mathrm{hyp}}(z,w)$, and from the 
$\mathrm{PSL}_{2}(\mathbb{R})$-invariance of $\hyp(z)$, formally for $w=u+iv\in X$, we find
\begin{align}\label{proofpropghypdelpinteqn1}
&\sum_{q\in \mathcal{P}}\sum_{\eta\in\Gamma_{q}\backslash \Gamma}\int_{X}g_{\mathrm{hyp}}(z,w)\del 
P_{\mathrm{gen},q}( \eta w)\hyp(w)=\notag\\&\sum_{q\in \mathcal{P}}\int_{0}^{\infty}\int_{0}^{1}g_{\mathrm{hyp}}
(z,\sigma_{q}w)\del P_{\mathrm{gen},q}(\sigma_{q}w)\frac{dudv}{v^{2}}.
\end{align}
Recall from equation \eqref{delP2}, that for any $w=u+iv\in\mathbb{H}$, the function 
$P_{\mathrm{gen},q}(\sigma_{q}w)$ does not depend on $u$. So the right-hand side of equation 
\eqref{proofpropghypdelpinteqn1} further decomposes to give 
\begin{align}
&\sum_{q\in \mathcal{P}}\int_{0}^{1\slash y_{p}}\int_{0}^{1}g_{\mathrm{hyp}}(z,\sigma_{q}w)
\del P_{\mathrm{gen},q}(\sigma_{q}w)\frac{dudv}{v^{2}}+\sum_{q\in \mathcal{P}}\int_{1\slash y_{p}}^{y_{p}}\bigg(\int_{0}^{1}
g_{\mathrm{hyp}}(z,\sigma_{q}w)du\bigg)\times\notag\\&\del P_{\mathrm{gen},q}(\sigma_{q}w)\frac{dv}{v^{2}}+\sum_{q\in \mathcal{P}}
\int_{y_{p}}^{\infty}\bigg(\int_{0}^{1}
g_{\mathrm{hyp}}(z,\sigma_{q}w)du\bigg)\del P_{\mathrm{gen},q}(\sigma_{q}w)\frac{dv}{v^{2}}.\label{decaux1}
\end{align} 
Since in the second line of formula (\ref{decaux1}) we have $1\slash y_{p}<v<y_{p}$, we can apply 
equation (\ref{zerocoeffecient1}), and rewrite the second line of formula (\ref{decaux1}) as 
\begin{align}
&\sum_{q\in \mathcal{P}}\int_{1\slash y_{p}}^{y_{p}}\bigg(4\pi\kappa_{q}(z)-\frac{4\pi}{\vx(X)}-
\frac{4\pi\log v}{\vx(X)}\bigg)\del P_{\mathrm{gen},q}(\sigma_{q}w)\frac{dv}{v^{2}}+\notag\\&4\pi\int_{1\slash y_{p}}
^{y_{p}}( v-y_{p})\del P_{\mathrm{gen},p}(\sigma_{p}w)\frac{dv}{v^{2}}.\label{propghypdelpinteqn1}
\end{align}
Since in the third line of formula (\ref{decaux1}) we have $v >y_{p}>1\slash y_{p}$, we can apply 
equation (\ref{zerocoeffecient2}), and rewrite the third line of formula (\ref{decaux1}) as 
\begin{align}
\sum_{q\in\mathcal{P} }\int_{y_{p}}^{\infty}\bigg(4\pi\kappa_{q}(z)-
\frac{4\pi}{\vx(X)}-\frac{4\pi\log v}{\vx(X)}\bigg)\del P_{\mathrm{gen},q}(\sigma_{q}w)\frac{dv}{v^{2}}.\label{propghypdelpinteqn2}
\end{align}
The proof of the proposition follows from combining equations (\ref{propghypdelpinteqn1}) and (\ref{propghypdelpinteqn2}).
\end{proof}
\end{prop}
\begin{rem}\label{remabsolute}
The formal unfolding of the integral obtained in Proposition \ref{propghypdelpint} translates into an equality of 
integrals, only if each of the three integrals on the right-hand side of equation (\ref{propghypdelpinteqn}) 
converges absolutely, which we prove in the lemmas that follow. 
\end{rem}
\begin{lem}\label{lem1}
Let $p,q\in \mathcal{P}$ be two cusps. For $z\in X$ and $w=u+iv\in\mathbb{H}$, the integral 
\begin{align*}
\int_{0}^{1\slash y_{p}}\int_{0}^{1}g_{\mathrm{hyp}}(z,\sigma_{q}w)\del P_{\mathrm{gen},q}
(\sigma_{q}w)\frac{dudv}{v^{2}} 
\end{align*}
converges absolutely. Furthermore as $z\in X$ approaches a cusp $p\in\mathcal{P}$, we have 
\begin{align}\label{lem1eqn}
\sum_{q\in \mathcal{P}}\int_{0}^{1\slash y_{p}}\int_{0}^{1}g_{\mathrm{hyp}}(z,\sigma_{q}w)\del P_{\mathrm{gen},q}(\sigma_{q}w)
\frac{dudv}{v^{2}} =o_{z}(1), 
\end{align}
where the contribution from the term $o_{z}(1)$ is a smooth function in $z$, which approaches zero, as 
$z\in X$ approaches the cusp $p\in\mathcal{P}$.
\begin{proof}
For $v\in \mathbb{R}_{> 0}$, from the formula for the function $\del P_{\mathrm{gen},q}(\sigma_{q}w)$ from 
equation \ref{delP2}, we derive that 
\begin{align*}
\frac{\del P_{\mathrm{gen},q}(\sigma_{q}w)}{v^{2}}=\frac{8\pi^{2}}{\sinh^{2}(2\pi v)}-\frac{2}{v^{2}}
\end{align*}
remains bounded. So it suffices to show that the integral 
\begin{align*}
\int_{0}^{1\slash y_{p}}\int_{0}^{1}g_{\mathrm{hyp}}(z,\sigma_{q}w)dudv
\end{align*}
converges absolutely. Let $\mathcal{I}$ denote the set $[0,1]\times [0,1\slash y_{p}]$. 
We view the above integral as a real-integral on the compact subset $\mathcal{I}\subset \mathbb{R}^{2}$. 
The hyperbolic Green's function $g_{\mathrm{hyp}}(z,\sigma_{q}w)$ is at most $\log$-singular on a measure zero 
subset of the interior points of $\mathcal{I}$. Furthermore from equation \eqref{ghyploglog}, the hyperbolic Green's 
function $g_{\mathrm{hyp}}(z,\sigma_{q}w)$ is at most $\log\log$-singular on a measure zero subset of the 
boundary points of $\mathcal{I}$. Hence, it is absolutely integrable on $\mathcal{I}$. This implies that the integral
\begin{align*}
\int_{0}^{1\slash y_{p}}\int_{0}^{1}g_{\mathrm{hyp}}(z,\sigma_{q}w)\del P_{\mathrm{gen},q}(\sigma_{q}w)
\frac{dudv}{v^{2}}
\end{align*}
converges absolutely, and also proves the asymptotic relation asserted in equation (\ref{lem1eqn}).
\end{proof}
\end{lem}
\begin{lem}\label{lem2}
Let $p,q\in \mathcal{P}$ be two cusps. For $z\in X$ and $w=u+iv\in\mathbb{H}$, the integral
\begin{align}\label{lem2int}
\int_{1\slash y_{p}}^{\infty}\bigg(4\pi\kappa_{q}(z)-\frac{4\pi}{\vx(X)}\bigg)\del P_{\mathrm{gen},q}(\sigma_{q}w) 
\frac{dv}{v^{2}} 
\end{align}
converges absolutely. Furthermore, as $z\in X$ approaches the cusp $p\in\mathcal{P}$, we have 
\begin{align*}
&\sum_{q\in  \mathcal{P}}\int_{1\slash y_{p}}^{\infty}\bigg(4\pi\kappa_{q}(z)-
\frac{4\pi}{\vx(X)}\bigg)\del P_{\mathrm{gen},q}(\sigma_{q}w)\frac{dv}{v^{2}}=\notag
\\&16\pi^{2}\bigg(-y_{p}+\frac{|\mathcal{P}|\big(\log y_{p}+1\big)}{\vx(X)}
-\sum_{q\in \mathcal{P}}k_{q,p}(0)+
\frac{2\pi}{3}\bigg)+O\bigg(\frac{\log y_{p}}{y_{p}}\bigg).
\end{align*}
\begin{proof}
Substituting the formula for the function $\del P_{\mathrm{gen},q}(\sigma_{q}w)$ from equation \eqref{delP1}, we have
\begin{align*}
&\int_{1\slash y_{p}}^{\infty}\bigg(4\pi\kappa_{q}(z)-\frac{4\pi}{\vx(X)}\bigg)\del P_{\mathrm{gen},q}(\sigma_{q}w) 
\frac{dv}{v^{2}} =\\&
\bigg(8\pi\kappa_{q}(z)-\frac{8\pi}{\vx(X)}\bigg)\int_{1\slash y_{p}}^{\infty}\bigg(\bigg(\frac{2\pi v }
{\sinh(2\pi v)}\bigg)^{2}-1\bigg)\frac{dv}{v^{2}}.
\end{align*}
The integral on the right-hand side of the above equation further simplifies to give
\begin{align}
&\bigg(8\pi\kappa_{q}(z)-\frac{8\pi}{\vx(X)}\bigg)\bigg{[}
\frac{1}{v}-2\pi\coth(2\pi v)\bigg{]}_{1\slash y_{p}}^{\infty}=\notag\\&
\bigg(8\pi\kappa_{q}(z)-\frac{8\pi}{\vx(X)}\bigg)\bigg(-2\pi-y_{p}+2\pi\coth\bigg(\frac{2\pi} {y_{p}}\bigg)\bigg).
\label{lem2eqn1}
\end{align}
Hence, from equation (\ref{lem2eqn1}), we can conclude that the integral (\ref{lem2int}) converges absolutely.

\vspace{0.2cm}
We now compute the asymptotics of the expression obtained on the right-hand side of equation (\ref{lem2eqn1}), as $z\in X$ 
approaches  the cusp $p\in \mathcal{P}.$ We first compute the asymptotics for the expression in the 
second bracket on the right-hand side of equation (\ref{lem2eqn1}).

For $t\in\mathbb{R}_{>0}$, recall that the Taylor series expansion of the function $\coth(t)$ as $t$ 
approaches zero is of the form
\begin{align*}
\coth(t)=\frac{1}{t}+\frac{t}{3}+O(t^3). 
\end{align*}
As $z\in X$ approaches $p\in \mathcal{P}$, the quantity $1\slash y_{p}$ approaches 
zero. So as $z\in X$ approaches $p\in \mathcal{P}$, using the Taylor expansion  of $\coth(2\pi\slash y_{p})$,  
we have the  asymptotic relation
\begin{align}
&-2\pi-y_{p}+ 2\pi\coth\bigg(\frac{2\pi}{y_{p}}\bigg)=
-2\pi-y_{p}+2\pi\bigg(\frac{y_{p}}{2\pi}+\frac{2\pi}{3y_{p}}+
O\bigg(\frac{1}{y_{p}^{3}}\bigg)\bigg)=-2\pi+\frac{4\pi^{2}}{3y_{p}}+
O\left(\frac{1}{y_{p}^{3}}\right).\label{lem2eqn2}
\end{align}
As $z\in X$ approaches $p\in \mathcal{P}$, from the Fourier expansion of Kronecker's limit function $\kappa_{q}(z)$ described in 
equation \eqref{fourierkappaeqn}, we have the following asymptotic relation 
\begin{align}\label{lem2eqn21}
&8\pi\kappa_{q}(z)-\frac{8\pi}{\vx(X)}=
8\pi \delta_{p,q} y_{p}-\frac{8\pi\log y_{p}}{\vx(X)}+8\pi k_{q,p}(0)-\frac{8\pi}{\vx(X)}
 +O\big(e^{-2\pi y_{p}}\big).
\end{align}
Combining equations (\ref{lem2eqn2}) and (\ref{lem2eqn21}), as $z\in X$ approaches 
$p\in \mathcal{P}$, we have the asymptotic relation 
for the right-hand side of equation (\ref{lem2eqn1})
\begin{align*}
\bigg(8\pi \delta_{p,q} y_{p}-\frac{8\pi\log y_{p}}{\vx(X)}+8\pi k_{q,p}(0)-\frac{8\pi}{\vx(X)}+
O\big(e^{-2\pi y_{p}}\big)\bigg)\bigg(-2\pi +\frac{4\pi^{2}}{3y_{p}}
+O\bigg(\frac{1}{y_{p}^{3}}\bigg)\bigg)=\\16\pi^{2}\bigg(-\delta_{p,q}y_{p}+
\frac{(\log y_{p}+1)}{\vx(X)}-k_{q,p}(0)+\frac{2\pi}{3}\delta_{p,q}+O\bigg(\frac{\log y_{p}}{y_{p}}\bigg)\bigg).
\end{align*}
Hence, taking the summation over all $q\in\mathcal{P}$ completes the proof of the lemma.
\end{proof}
\end{lem}
\begin{lem}\label{lem3}
Let $p,q\in \mathcal{P}$ be two cusps. For $z\in X$ and $w=u+iv\in\mathbb{H}$, the integral 
\begin{align}\label{lem3inteqn1}
\int_{1\slash y_{p}}^{\infty}\frac{4\pi\log v}{\vx(X)}\del P_{\mathrm{gen},q}(\sigma_{q}w)\frac{dv}{v^{2}}
\end{align}
converges absolutely. Furthermore, we have the upper bound
\begin{align}\label{lem3inteqn2}
\sum_{q\in \mathcal{P}}\int_{1\slash y_{p}}^{\infty}
\Bigg{|}\frac{4\pi\log v}{\vx(X)}\del 
P_{\mathrm{gen},q}(\sigma_{q}w)\Bigg{|}\frac{dv}{v^{2}}\leq \frac{8\pi|\mathcal{P}|}{\vx(X)}
\bigg(1+ \frac{4\pi^{2}}{3}\bigg).
\end{align}
\begin{proof}
We prove the upper bound asserted in (\ref{lem3inteqn2}), which also proves the absolute convergence 
of the integral in (\ref{lem3inteqn1}). Observing the elementary estimate
\begin{align}
&\int_{1\slash y_{p}}^{\infty}\bigg{|}\log v\del P_{\mathrm{gen},q}
(\sigma_{q}w)\bigg{|}\frac{dv}{v^{2}}\leq\notag\\&
\int_{0}^{1}\bigg{|}\log v\del P_{\mathrm{gen},q}(\sigma_{q}w)\bigg{|}\frac{dv}{v^{2}}+
\int_{1}^{\infty}\bigg{|}\log v\del P_{\mathrm{gen},q}(\sigma_{q}w)\bigg{|}\frac{dv}{v^{2}},\label{lem3inteqn10}
\end{align}
we proceed to bound the two integrals on the right-hand side of the above inequality. For $v\in\mathbb{R}_{> 0}$, from the 
equation \eqref{delP1}, we find that the function 
\begin{align*}
-\frac{\del P_{\mathrm{gen},q}(\sigma_{q}w)}{v^{2}}=\frac{2}{v^{2}}-\frac{8\pi^{2} }{\sinh^{2}(2\pi  v)}
\end{align*}
is a positive monotone decreasing function, and hence, attains its maximum value at $v= 0$. So we compute the limit  
\begin{align*}
-\lim_{v\rightarrow 0}\frac{\del P_{\mathrm{gen},q}(\sigma_{q}w)}{v^{2}}
=\lim_{v\rightarrow 0}\bigg(\frac{2}{v^{2}}-\frac{8\pi^{2} }{\sinh^{2}(2\pi  v)}\bigg)=\frac{8\pi^{2}}{3}.
\end{align*}
So using the fact that, for $v\in (0,1]$, $|\log v|=-\log v$, we have the following upper bound for the first integral on the 
right-hand side of inequality (\ref{lem3inteqn10})
\begin{align}
&\int_{0}^{1}\bigg{|}\log v\del P_{\mathrm{gen},q}(\sigma_{q}w)\bigg{|}\frac{dv}{v^{2}}\leq 
-\frac{8\pi^{2}}{3}\int_{0}^{1}\log v dv=\frac{8\pi^{2}}{3}.\label{lem3eqn2}
\end{align}
Again using formula \eqref{delP1}, we derive
\begin{align*}
 \max_{v\in\mathbb{R}_{>0}}\big{|}\del P_{\mathrm{gen},q}(\sigma_{q}w)\big{|}=
 \max_{v\in\mathbb{R}_{>0}}\bigg(2-\frac{8\pi^{2}v^{2} }{\sinh^{2}(2\pi  v)} \bigg)=2.
\end{align*}
Using the above bound, we derive the following upper bound for the second integral on the right-hand side of 
inequality (\ref{lem3inteqn10})
\begin{align}
&\int_{1}^{\infty}\bigg{|}\log v\del P_{\mathrm{gen},q}(\sigma_{q}w)\bigg{|}\frac{dv}{v^{2}}\leq2\int_{1}^{\infty}\frac{\log v }{v^{2}}dv=2\bigg(\bigg{[}-\frac{\log v}{v}\bigg{]}_{1}^{\infty}+
\bigg{[}-\frac{1}{v}\bigg{]}_{1}^{\infty}\bigg)=2.\label{lem3eqn3}
\end{align}
Hence, combining the upper bounds derived in equations (\ref{lem3eqn2}) and (\ref{lem3eqn3}) proves the lemma.
\end{proof}
\end{lem}
\begin{lem}\label{lem4}
Let $p\in\mathcal{P}$ be a cusp. For $z\in X$ and $w=u+iv\in\mathbb{H}$ with $y_{p}> 1$, 
the integral  
\begin{align}\label{lem4int}
-4\pi y_{p}\int_{1\slash y_{p}}^{y_{p}}\del P_{\mathrm{gen},p}(\sigma_{p}w)\frac{dv}{v^{2}}
\end{align}
converges absolutely. Furthermore, as $z\in X$ approaches the cusp $p\in \mathcal{P}$, we have 
\begin{align}
-&4\pi y_{p}\int_{1\slash y_{p}}^{y_{p}}\del P_{\mathrm{gen},p}(\sigma_{p}w)\frac{dv}{v^{2}}=
4\pi\bigg(4\pi y_{p}\coth(2\pi y_{p})-2-\frac{8\pi^{2}}{3}\bigg) + O\bigg(\frac{1}{y_{p}^{2}}
\bigg). \label{lem4eqn}
\end{align}
\begin{proof}
From equation \eqref{delP1}, for a cusp $p\in \mathcal{P}$, we find 
\begin{align}\label{lem4eqn2}
-&4\pi y_{p}\int_{1\slash y_{p}}^{y_{p}}\del P_{\mathrm{gen},p}(\sigma_{p}w)
\frac{dv}{v^{2}}= -8\pi y_{p}\int_{1\slash y_{p}}^{y_{p}}\bigg(\frac{4\pi^2}{\sinh^{2}(2\pi v)}-\frac{1}{v^2}\bigg)dv
=\notag\\-&8\pi y_{p}\bigg{[}\frac{1}{v}-2\pi\coth(2\pi v)\bigg{]}_{1\slash y_{p}}^{y_{p}}=-
8\pi y_{p}\bigg(\frac{1}{y_{p}}-2\pi\coth(2\pi y_{p})-y_{p}+2\pi\coth\bigg(\frac{2\pi}{y_{p}}\bigg)\bigg)
=\notag\\-&8\pi +16\pi^{2}y_{p}\coth(2\pi y_{p})-8\pi y_{p}\bigg(-y_{p}+ 2\pi\coth\bigg(\frac{2\pi}{y_{p}}\bigg) \bigg).
\end{align}
This implies that the integral (\ref{lem4int}) converges absolutely.

As $z\in X$ approaches the cusp $p\in \mathcal{P}$, from the Taylor expansion of 
$\coth(2\pi\slash y_{p})$ already used in equation (\ref{lem2eqn2}), we get
\begin{align*}
-8\pi y_{p}\bigg(-y_{p}+ 2\pi\coth\bigg(\frac{2\pi}{y_{p}}\bigg)\bigg)=-8\pi y_{p}\bigg(
\frac{4\pi^{2}}{3 y_{p}}+O\bigg(\frac{1}{y_{p}^{3}}\bigg)\bigg)
-\frac{32\pi^{3}}{3 }+O\bigg(\frac{1}{y_{p}^{2}}\bigg),
\end{align*}
which together with equation \eqref{lem4eqn2} completes the proof of the lemma.
\end{proof}
\end{lem}
\begin{lem}\label{lem5}
Let $p\in\mathcal{P}$ be a cusp. For $z\in X$ and $w=u+iv\in\mathbb{H}$ with $y_{p}> 1$, the integral   
\begin{align}\label{lem5int}
4\pi\int_{1\slash y_{p}}^{y_{p}} \del P_{\mathrm{gen},p}(\sigma_{p}w)\frac{dv}{v}. 
\end{align}
converges absolutely. Furthermore, as $z\in X$ approaches the cusp $p\in\mathcal{P}$, we have
\begin{align}\label{lem5eqn}
4\pi\int_{1\slash y_{p}}^{y_{p}} \del P_{\mathrm{gen},p}(\sigma_{p}w)\frac{dv}{v}= 
-8\pi\log y_{p}+8\pi(1-\log (4\pi))+O\bigg(\frac{1}{y_{p}}\bigg).
\end{align}
\begin{proof}
Using equation \eqref{delP1}, for a cusp $p\in \mathcal{P}$, we find 
\begin{align*}
&\int_{1\slash y_{p}}^{y_{p}}\del P_{\mathrm{gen},p}(\sigma_{p}w)\frac{dv}{v}=
2\int_{1\slash y_{p}}^{y_{p}}\bigg(-\frac{1}{v}+\frac{4\pi^{2}v}{\sinh^{2}(2\pi v)}\bigg)dv=\\&
2\bigg{[}-\log v-2\pi v\coth(2\pi v)+\log(\sinh(2\pi v))
\bigg{]}_{1\slash y_{p}}^{y_{p}}.
\end{align*}
Substituting the formulae for $\coth(2\pi v)$ and $\sinh(2\pi v)$, the right-hand side of the above 
equation can be further simplified to 
\begin{align*}
2\bigg{[}-\log v-4\pi v-\frac{4\pi v}{e^{4\pi v}-1}+\log\bigg(\frac{e^{4\pi v}-1}{2}\bigg)\bigg{]}_
{1\slash y_{p}}^{y_{p}}.
\end{align*}
Observe that  
\begin{align}
&\bigg{[}-\log v-4\pi v-\frac{4\pi v}{e^{4\pi v}-1}+\log\bigg(\frac{e^{4\pi v}-1}{2}\bigg)\bigg{]}_
{1\slash y_{p}}^{y_{p}}=-\log y_{p}-4\pi y_{p}-\frac{4\pi y_{p}}{e^{4\pi y_{p}}-1}+\notag\\&
\log\big(e^{4\pi y_{p}}-1\big)+\log\bigg(\frac{1}{y_{p}}\bigg)+ \frac{4\pi}{y_{p}}+
\frac{4\pi}{y_{p}\big(e^{4\pi\slash y_{p}}-1\big)}-\log\big(e^{4\pi\slash y_{p}}-1\big)=\notag\\
-&\log y_{p}-\log\bigg(\frac{e^{4\pi y_{p}}}{e^{4\pi y_{p}}-1}\bigg)-\frac{4\pi y_{p}}{e^{4\pi y_{p}}-1}+
\frac{4\pi}{y_{p}}+\frac{4\pi}{y_{p}\big(e^{4\pi\slash y_{p}}-1\big)}-\log\big(y_{p}\big(e^{4\pi\slash y_{p} }-
1\big)\big) ,\label{lem5eqn1}
\end{align}
which proves that the integral (\ref{lem5int}) converges absolutely. 

We now compute the asymptotic expansion of each of the terms in the above 
expression, as $z\in X$ approaches the cusp $p\in\mathcal{P}$. As $z\in X$ approaches 
the cusp $p\in\mathcal{P}$, we have the asymptotic relation for the 
first and second terms of (\ref{lem5eqn1})
\begin{align}\label{lem5eqn2}
&-\log y_{p}-\log\bigg(\frac{e^{4\pi y_{p}}}{e^{4\pi y_{p}}-1}\bigg)=-\log y_{p}+O\left(e^{-4\pi y_{p}}\right);
\end{align}
the third and fourth terms of (\ref{lem5eqn1}) satisfy the asymptotic relation 
\begin{align}\label{lem5eqn3}
-\frac{4\pi y_{p}}{e^{4\pi y_{p}}-1}+\frac{4\pi}{y_{p}}=O\bigg(\frac{1}{ y_{p}}\bigg);
\end{align}
the fifth term satisfies the asymptotic relation 
\begin{align}\label{lem5eqn5}
&\frac{4\pi}{y_{p}
\big(e^{4\pi\slash y_{p}}-1\big)}=\frac{4\pi}{y_{p}\bigg(\displaystyle\sum_{n=1}^{\infty}
\frac{(4\pi)^{n}}{n!\,y_{p}^{n}}\bigg)}=1+O\bigg(\frac{1}{y_{p}}\bigg);
\end{align}
and the sixth term satisfies the asymptotic relation 
\begin{align}\label{lem5eqn4}
&-\log\big(y_{p}\big(e^{4\pi\slash y_{p}} -1\big)\big)=-\log\bigg(
\sum_{n=1}^{\infty}\frac{(4\pi)^{n}}{n!\,y_{p}^{n-1}}\bigg)=\notag\\&
-\log\bigg(4\pi+\sum_{n=1}^{\infty}\frac{(4\pi)^{n+1}}{(n+1)!\,y_{p}^{n}}\bigg)=
-\log( 4\pi) + O\bigg(\frac{1}{y_{p}}\bigg).
\end{align} 
Substituting the asymptotic relations obtained in equations (\ref{lem5eqn2}), (\ref{lem5eqn3}), 
(\ref{lem5eqn5}), and (\ref{lem5eqn4}) into (\ref{lem5eqn1}), we derive the asymptotic relation 
\begin{align*}
4\pi\int_{1\slash y_{p}}^{y_{p}} \del P_{\mathrm{gen},p}(\sigma_{p}w)\frac{dv}{v}= 
-8\pi\log y_{p}+8\pi(1-\log (4\pi))+O\bigg(\frac{1}{y_{p}}\bigg),
\end{align*}
as $z\in X$ approaches the cusp $p\in \mathcal{P}$, which completes the proof of the lemma.
\end{proof}
\end{lem}
In the following proposition, combining all the asymptotics established in this section, 
we compute the asymptotics of the integral
\begin{align*}
\int_{X}g_{\mathrm{hyp}}(z,w)\del P(w)\hyp(w), 
\end{align*}
as $z\in X$ approaches a cusp $p\in\mathcal{ P}.$
\begin{prop}\label{prop2}
Let $p\in\mathcal{P}$ be a cusp. Then, as $z\in X$ approaches the cusp $p\in\mathcal{P}$, we have
\begin{align}
&\int_{X}g_{\mathrm{hyp}}(z,w)\del P(w)\hyp(w)=-\notag\\&\frac{32\pi^{2}(g-1)\log y_{p}}{\vx(X)}-
\sum_{q\in \mathcal{P}}\int_{1\slash y_{p}}^{\infty}\frac{4\pi\log v}{\vx(X)}\del P_{\mathrm{gen},q}
(\sigma_{q}w)\frac{dv}{v^{2}}+\alpha_{p}+o_{z}(1),\notag\\&\mathrm{where}\,\,\alpha_{p}=\frac{16\pi^{2}|\mathcal{P}|}{
\vx(X)}-16\pi^{2}\sum_{q\in\mathcal{P}}k_{q,p}(0)-8\pi\log (4\pi),   \label{alphapdefn}
\end{align}
and the contribution from the term $o_{z}(1)$ is a smooth function in $z$, which approaches zero, as 
$z\in X$ approaches the cusp $p\in\mathcal{P}$.
\begin{proof}
From Lemmas \ref{lem1}, \ref{lem2}, \ref{lem3}, \ref{lem4}, and \ref{lem5}, it follows that each of the integrals on 
the right-hand side of the equation \eqref{propghypdelpinteqn} is absolutely convergent. This implies that the 
equality of integrals described in equation \eqref{propghypdelpinteqn} indeed holds true for all $z\in X$ provided 
that $y_{p}>1$. 

As $z\in X$ approaches the cup $p\in \mathcal{P}$, combining Lemmas 
\ref{lem1} and \ref{lem2}, we find that the first two integrals on the 
right-hand side of equation \eqref{propghypdelpinteqn} yield
\begin{align}
&16\pi^{2}\bigg(-y_{p}+\frac{|\mathcal{P}|\big(\log y_{p}+1\big)}{\vx(X)}
-\sum_{q\in \mathcal{P}}k_{q,p}(0)+
\frac{2\pi}{3}\bigg)-\notag\\&\sum_{q\in \mathcal{P}}
\int_{1\slash y_{p}}^{\infty}\frac{4\pi\log v}{\vx(X)}\del P_{\mathrm{gen},q}(\sigma_{q}w)
\frac{dv}{v^{2}}+ o_{z}(1),\label{prop2eqn2}
\end{align}
where the contribution from the term $o_{z}(1)$ is a smooth function in $z$, which approaches zero, as 
$z\in X$ approaches $p\in\mathcal{P}$. As $z\in X$ approaches the cusp $p\in \mathcal{P}$, combining Lemmas 
\ref{lem4} and \ref{lem5}, we find that the third integral on the right-hand side of 
equation \eqref{propghypdelpinteqn} yields
\begin{align}\label{prop2eqn3}
16\pi^{2}y_{p}\coth(2\pi y_{p})-8\pi\log y_{p}-\frac{32\pi^{3}}{3}-8\pi\log (4\pi)+
O\bigg(\frac{1}{y_{p}}\bigg).
\end{align}
Combining (\ref{prop2eqn2}) and (\ref{prop2eqn3}), as $z\in X$ approaches the 
cusp $p\in \mathcal{P}$, the right-hand side of equation \eqref{propghypdelpinteqn} simplifies to 
\begin{align}
-&16\pi^{2}y_{p}+16\pi^{2}y_{p}\coth(2\pi y_{p})+\frac{16\pi^{2}|\mathcal{P}|
\log y_{p}}{\vx(X)}-8\pi\log y_{p}+\frac{16\pi^{2}|\mathcal{P}|}{\vx(X)}-
\notag\\&16\pi^{2}\sum_{q\in \mathcal{P}}k_{q,p}(0)-\sum_{q\in \mathcal{P}}
\int_{1\slash y_{p}}^{\infty}\frac{4\pi\log v}{\vx(X)}\del P_{\mathrm{gen},q}(\sigma_{q}w)
\frac{dv}{v^{2}}-8\pi\log (4\pi)+ o_{z}(1).\label{prop2eqn4}
\end{align}
As $z\in X$ approaches the cusp $p\in\mathcal{P}$, we have the following asymptotic relation for the first two terms 
in the above expression
\begin{align*}
&16\pi^{2}y_{p}\big(\coth(2\pi y_{p})-1\big)=16\pi^{2}y_{p}
\bigg(\frac{\cosh(2\pi y_{p})-\sinh(2\pi y_{p})}{\sinh(2\pi y_{p})}\bigg)=
O\big(e^{-y_{p}}\big).
\end{align*}
Furthermore, as $z\in X$ approaches the cusp $p\in\mathcal{P}$ the third and fourth terms in expression 
\eqref{prop2eqn4} give
\begin{align*}
 \frac{16\pi^{2}|\mathcal{P}|\log y_{p}}{\vx(X)}-8\pi\log y_{p}=
 -\frac{32\pi^{2}(g-1)\log y_{p}}{\vx(X)}.
\end{align*}
Hence, as $z\in X$ approaches the cusp $p\in\mathcal{P}$, 
the expression in (\ref{prop2eqn4}) further reduces to give
\begin{align*}
-&\frac{32\pi^{2}(g-1)\log y_{p}}{\vx(X)} -\sum_{q\in \mathcal{P}}
\int_{1\slash y_{p}}^{\infty}\frac{4\pi\log v}{\vx(X)}\del P_{\mathrm{gen},q}(\sigma_{q}w)
\frac{dv}{v^{2}}+\notag\\&\frac{16\pi^{2}|\mathcal{P}|}{\vx(X)}-16\pi^{2}\sum_{q\in\mathcal{P}}k_{q,p}(0)-
8\pi\log (4\pi)+ o_{z}(1),
\end{align*}
which completes the proof of the proposition.
\end{proof}
\end{prop}
\begin{cor}\label{cor1}
Let $p\in \mathcal{P}$ be a cusp. Then, as $z\in X$ 
approaches the cusp $p\in\mathcal{P}$, we have
\begin{align*}
\phi(z)=-&\frac{4\pi\log y_{p}}{\vx(X)}-\sum_{q\in \mathcal{ P}}
\int_{1\slash y_{p}}^{\infty}\frac{\log v}{2g\vx(X)}\del P_{\mathrm{gen},q}(\sigma_{q}w)
\frac{dv}{v^{2}}+\notag\\&\frac{\alpha_{p}}{8\pi g}+\frac{2\pi k_{p,p}(0)}{g}-
\frac{C_{\mathrm{hyp}}}{8g^{2}}-\frac{2\pi c_{X}}{g\vx(X)}+o_{z}(1),
\end{align*}
where the constant $\alpha_{p}$ is as defined in (\ref{alphapdefn}), and the contribution from the term $o_{z}(1)$ is a 
smooth function in $z$, which approaches zero, as $z\in X$ approaches the cusp $p\in\mathcal{P}$. 
\begin{proof}
As $z\in X$ approaches the cusp $p\in\mathcal{P}$, from equation \eqref{H(z)cusp}, we have
\begin{align}
&\frac{H(z)}{2g}= -\frac{4\pi \log y_{p}}{g\vx(X)}-\frac{2\pi}{g\vx(X)}+\frac{2\pi
k_{p,p}(0)}{g}+O\bigg(\frac{1}{y_{p}}\bigg).\label{cor1eqn1}
\end{align}
Furthermore, from Proposition \ref{prop2}, we find that
\begin{align}
&\frac{1}{8\pi g}\int_{X}g_{\mathrm{hyp}}(z,\zeta)\del P(\zeta)\hyp(\zeta)=
-\frac{4\pi\log y_{p}}{\vx(X)}+\frac{4\pi \log y_{p}}{g\vx(X)}-\notag\\&\sum_{q\in
\mathcal{P}}\int_{1\slash y_{p}}^{\infty}
\frac{\log v}{2g\vx(X)}\del P_{\mathrm{gen},q}(\sigma_{q}w)\frac{dv}{v^{2}}+
\frac{\alpha_{p}}{8\pi g} +o_{z}(1),\label{cor1eqn2}
\end{align}
where the contribution from the term $o_{z}(1)$ is a smooth function in $z$, which approaches zero, as 
$z\in X$ approaches the cusp $p\in\mathcal{P}$. The proof of the corollary follows from 
combining equations \eqref{phi(z)formula}, (\ref{cor1eqn1}), and (\ref{cor1eqn2}).
\end{proof}
\end{cor}
The following proposition has been proved as Proposition 6.1.9 in \cite{anilthesis} (or Proposition 4.10 in \cite{anilpaper2}). 
However, for the convenience of the reader, we reproduce the proof here. 
\begin{prop}\label{prop3}
We have the following upper bound
\begin{align*}
\frac{\big|C_{\mathrm{hyp}}\big|}{8g^{2}} \leq \frac{2\pi \left(d_{X}+1\right)^{2}}{\lambda_{1}
\vx(X)},
\end{align*}
where $\lambda_{1}$ denotes the first non-zero eigenvalue of the hyperbolic Laplacian acting on smooth functions 
defined on $X$. 
\begin{proof}
Recall that $C_{\mathrm{hyp}}$ is defined as 
\begin{align*}
&C_{\mathrm{hyp}}= \int_{X}\int_{X}\ghyp(\zeta,\xi)\bigg(\int_{0}^{\infty}\del 
\khyp(t;\zeta)dt\bigg)\bigg(\int_{0}^{\infty}\del \khyp(t;\xi)dt\bigg)\hyp(\xi)\hyp(\zeta).
\end{align*}
From formulae \eqref{phi}, \eqref{phi(z)formula}, we have
\begin{align}
&\del\phi(z)=\frac{4\pi\can(z)}{\hyp(z)}-\frac{4\pi}{\vx(X)}\Longrightarrow \int_{X}\del\phi(z)
\hyp(z)=0,\label{prop3.26eqn0}\\&\phi(z)= \frac{1}{2g}\int_{X}g_{\mathrm{hyp}}(z,\zeta)\bigg(
\int_{0}^{\infty}\del K_{\mathrm{hyp}}(t;\zeta)dt\bigg)\hyp(\zeta)-\frac{C_{\mathrm{hyp}}}{8g^{2}},\notag
\end{align}
respectively. So combining the above two equations, we get
\begin{align}
-&\frac{1}{4\pi}\int_{X}\phi(z)\del\phi(z)\hyp(z)=\frac{1}{2g}\int_{X}\int_{X}
g_{\mathrm{hyp}}(z,\zeta)\bigg(\int_{0}^{\infty}\del K_{\mathrm{hyp}}(t;\zeta)dt\bigg)\hyp(\zeta)
\can(z).\label{prop3.26eqn2}
\end{align}
Observe that
\begin{align*}
\int_{X}g_{\mathrm{hyp}}(z,\zeta)\bigg(\int_{0}^{\infty}\del K_{\mathrm{hyp}}(t;\zeta)dt\bigg)
\hyp(\zeta)=2g\phi(z)+\frac{C_{\mathrm{hyp}}}{4g}\in C_{\ell,\ell\ell}(X).
\end{align*}
So using equations \eqref{keyidentity} and \eqref{prop3.26eqn2}, we derive
\begin{align}
&\int_{X}\phi(z)\del\phi(z)\hyp(z)=\frac{\pi}{ g^{2}}\int_{X}\int_{X}
g_{\mathrm{hyp}}(z,\zeta)\left(\int_{0}^{\infty}\del K_{\mathrm{hyp}}(t;\zeta)dt\right)\times\notag
\\&\left(\int_{0}^{\infty}\del K_{\mathrm{hyp}}(t;z)dt\right)
\hyp(\zeta)\hyp(z)=\frac{\pi C_{\mathrm{hyp}}}{ g^{2}}.\label{prop3.26eqn1}
\end{align}
From equation \eqref{prop3.26eqn0}, we have
\begin{align}\label{prop3.26eqn3}
&\sup_{z\in X}|\del\phi(z)|\leq \sup_{z\in X}\bigg{|}\frac{4\pi\can(z)}{\vx(X)\shyp(z)}\bigg{|} +
\frac{4\pi}{\vx(X)}=\frac{4\pi\left(d_{X}+1\right)}{\vx(X)},
\end{align}
where $d_{X}$ is as defined in \eqref{defndx}. As the function $\phi(z)\in L^{2}(X)$, it admits a 
spectral expansion in terms of the eigenfucntions of the hyperbolic Laplacian $\del$. So from the arguments used to prove Proposition 4.1 in \cite{jkannals}, we have 
\begin{align}\label{prop3.26eqn4}
&\bigg|\int_{X}\phi(z)\del \phi(z)\hyp(z) \bigg|
\leq\sup_{z\in X} \frac{|\del\phi(z)|^{2}}{\lambda_{1}}\int_{X}\hyp(z),
\end{align}
where $\lambda_{1}$ denotes the first non-zero eigenvalue of the hyperbolic Laplacian $\del$. Hence, from equation 
\eqref{prop3.26eqn1}, and combining estimates \eqref{prop3.26eqn3} and 
\eqref{prop3.26eqn4}, we arrive at the estimate
\begin{align*}
\big|C_{\mathrm{hyp}}\big| =\frac{ g^{2}}{\pi}\bigg| \int_{X}\phi(z)\del \phi(z) \hyp(z)
\bigg|\leq\frac{g^{2}}{\pi\lambda_{1}}\int_{X}|\del \phi(z)|^{2}\hyp(z)\leq\frac{16\pi 
g^{2}\left(d_{X}+1\right)^{2}}{\lambda_{1}\vx(X)},
\end{align*}
which completes the proof of the proposition.
\end{proof}
\end{prop}
\begin{thm}\label{finalthm}
Let $p,q \in\mathcal{ P}$ be two cusps with $p\not = q$. Then, we have the upper bound
\begin{align*}
&\big{|} g_{\mathrm{can}}(p,q)\big{|}\leq 4\pi \big{|}k_{p,q}(0)\big{|}+
\frac{2\pi}{g}\bigg(\sum_{\substack{s\in \mathcal{P}\\s\not = p}}\big{|}k_{s,p}(0)\big{|}+
\sum_{\substack{s\in  \mathcal{P}\\s\not = q}}\big{|}k_{s,q}(0)\big{|}\bigg)+\notag\\&
\frac{1}{\vx(X)}\bigg(\frac{4\pi(d_{X}+1)^{2}}{\lambda_{1}}+\frac{\big{|}4\pi c_{X}\big{|}}{g}+
\frac{ 43|\mathcal{P}|}{g}+4\pi\bigg)+\frac{2\log (4\pi)}{g}.
\end{align*}
\begin{proof}
For $z,w\in X$, from equation \eqref{phi}, we have
\begin{align*}
g_{\mathrm{can}}(p,q)=\lim_{z\rightarrow p}\lim_{w\rightarrow q}\big(
g_{\mathrm{hyp}}(z,w)-\phi(z)-\phi(w)\big).
\end{align*}
Combining equation \eqref{ghypcusp} with Corollary \ref{cor1}, for a fixed $w\in X$ with  $z\in X$ approaching the 
cusp $p\in \mathcal{P}$, we have
\begin{align}
&\lim_{z\rightarrow p}\big(g_{\mathrm{hyp}}(z,w)-\phi(z)\big)=4\pi\kappa_{p}(w)-
\frac{4\pi}{\vx(X)}-\frac{\alpha_{p}}{8\pi g}-\frac{2\pi k_{p,p}(0)}{g}+\notag\\&
\frac{C_{\mathrm{hyp}}}{8g^{2}}+\frac{2\pi c_{X}}{g\vx(X)}+
\lim_{y_{p}\rightarrow\infty }\sum_{s\in\mathcal{P}}\int_{1\slash y_{p}}^{\infty}
\frac{\log\zeta}{2g\vx(X)}\del P_{\mathrm{gen},s}(\sigma_{s}\xi)\frac{d\zeta}{\zeta^{2}},\label{finalthmeqn1}
\end{align}
where $\zeta=\Im(\xi).$ As $w\in X$ approaches the cusp $q\in\mathcal{ P}$ with $q\not = p$, from the Fourier expansion of the Kronecker's limit function 
$\kappa_{p}(w)$, stated in equation \ref{fourierkappaeqn}, we have
\begin{align*}
4\pi\kappa_{p}(w)=4\pi k_{p,q}(0) -\frac{4\pi \log v_{q}}{\vx(X)}+O\big(e^{-2\pi v_{q}}
\big).
\end{align*}
So using Corollary \ref{cor1} one more time, and substituting the above asymptotic 
relation into equation (\ref{finalthmeqn1}), we compute the limit
\begin{align}
&\lim_{z\rightarrow p}\lim_{w\rightarrow q}\big(g_{\mathrm{hyp}}(z,w)-\phi(z)-\phi(w)
\big)= 4\pi k_{p,q}(0)-\frac{4\pi}{\vx(X)}-\frac{\alpha_{p}}{8\pi g}-\frac{2\pi k_{p,p}(0)}{g}
-\notag\\&\frac{\alpha_{q}}{8\pi g}-\frac{2\pi k_{q,q}(0)}{g}+\frac{C_{\mathrm{hyp}}}{4g^{2}}
+\frac{4\pi c_{X}}{g\vx(X)}+\lim_{y_{p}\rightarrow\infty }\sum_{s\in\mathcal{P}}\int_{1\slash y_{p}}^{\infty}
\frac{\log\zeta}{2g\vx(X)}\del P_{\mathrm{gen},s}(\sigma_{s}\xi)\frac{d\zeta}{\zeta^{2}}+
\notag\\&\hspace{6.4cm}\lim_{v_{q}\rightarrow\infty }\sum_{s\in\mathcal{ P}}
\int_{1\slash v_{q}}^{\infty}\frac{\log\zeta}{2g\vx(X)}
\del P_{\mathrm{gen},s}(\sigma_{s}\xi)\frac{d\zeta}{\zeta^{2}}.
\label{finalthmeqn2}
\end{align}
Using the definition of the constant $\alpha_{p}$ from (\ref{alphapdefn}), we find that 
the first six terms on the right-hand side of the above equation give
\begin{align*}
&4\pi k_{p,q}(0)-\frac{1}{ g}\bigg(\frac{2\pi|\mathcal{P}|}{\vx(X)}-2\pi\sum_{s\in \mathcal{P}}
k_{s,p}(0)-\log (4\pi)\bigg)-\frac{2\pi k_{p,p}(0)}{g}-\notag
\\&\frac{4\pi}{\vx(X)}-\frac{1}{g}\bigg(\frac{2\pi|\mathcal{P}|}{\vx(X)}-2\pi
\sum_{s\in\mathcal{P}}k_{s,q}(0)-\log (4\pi)\bigg)-\frac{2\pi k_{q,q}(0)}{g}=\notag\\
&4\pi k_{p,q}(0)+\frac{2\pi}{g}\bigg(\sum_{\substack{s\in \mathcal{P}\\s\not = p}}
k_{s,p}(0)+\sum_{\substack{s\in \mathcal{P}\\s\not = q}}k_{s,q}(0)\bigg)
-\frac{4\pi(|\mathcal{P}|+g)}{g\vx(X)}+\frac{2\log (4\pi)}{g}.
\end{align*}
Furthermore, the expression on the right-hand side of the above equation can be bounded by
\begin{align}
4\pi\big{|}k_{p,q}(0)\big{|}+\frac{2\pi}{g}\bigg(\sum_{\substack{s\in\mathcal{ P}
\\s\not = p}}\big{|}k_{s,p}(0)\big{|}+\sum_{\substack{s\in \mathcal{P}\\s\not = q}}\big{|}
k_{s,q}(0)\big{|}\bigg)+\frac{13|\mathcal{P}|+4\pi g}{g\vx(X)}+
\frac{2\log (4\pi)}{g}.\label{finalthmeqn3}
\end{align}
Using Proposition \ref{prop3}, we derive the upper bound for the next two terms on the right-hand 
side of equation (\ref{finalthmeqn2})
\begin{align}\label{finalthmeqn4}
\frac{C_{\mathrm{hyp}}}{4g^{2}}+\frac{4\pi c_{X}}{g\vx(X)} \leq 
\frac{4\pi\left(d_{X}+1\right)^{2}}{\lambda_{1}\vx(X)}+\frac{\big{|}4\pi c_{X}\big{|}}{g\vx(X)} .
\end{align}
From Lemma \ref{lem3}, we have the upper bound for the absolute value of the 
last two terms on the right-hand side of equation (\ref{finalthmeqn2})
\begin{align}
&\lim_{y_{p}\rightarrow\infty }\sum_{s\in\mathcal{P}}\int_{1\slash y_{p}}^{\infty}
\Bigg{|}\frac{\log\zeta}{2g\vx(X)}\del P_{\mathrm{gen},s}(\sigma_{s}\xi)\Bigg{|}\frac{d\zeta}{\zeta^{2}}
+\notag\\&\lim_{v_{q}\rightarrow\infty }\sum_{s\in \mathcal{P}}\int_{1\slash v_{q}}^{\infty}
\Bigg{|}\frac{\log\zeta}{2g\vx(X)}\del P_{\mathrm{gen},s}(\sigma_{s}\xi)
\Bigg{|}\frac{d\zeta}{\zeta^{2}}\leq 
\frac{2|\mathcal{P}|}{g\vx(X)}\bigg(1+ \frac{4\pi^{2}}{3}\bigg)\leq \frac{30|\mathcal{P}|}{g\vx(X)}.\label{finalthmeqn5}
\end{align}
The proof of the theorem follows from combining the estimates obtained in equations (\ref{finalthmeqn3}), 
(\ref{finalthmeqn4}), and (\ref{finalthmeqn5}).
\end{proof}
\end{thm}

{\small{}}
\vspace{0.3cm}
{\small{
Department of Mathematics, \\University of Hyderabad, \\Prof. C.~R.~Rao Road, Gachibowli,\\
Hyderabad, 500046, India\\email: anilatmaja@gmail.com}}
\end{document}